\documentclass[11pt]{article}

\usepackage{amsmath,amsthm,amssymb}
\usepackage{graphicx,color}
\usepackage{setspace}
\usepackage{hyperref}

\newcommand{\ra}{\rightarrow}
\newcommand{\R}{\mathbb{R}}

\newcommand{\abs}[1]{\lvert#1\rvert}
\newcommand{\norm}[1]{\lVert#1\rVert}

\makeatletter
\newcommand{\ps@paper}{%
  \renewcommand{\@oddhead}{%
    {\itshape Sum of squares law}\hfill \thepage}%
  \renewcommand{\@evenhead}{%
    {\itshape Sum of squares law}\hfill \thepage}%
  \renewcommand{\@oddfoot}{}%
  \renewcommand{\@evenfoot}{}}
\makeatother

\theoremstyle{plain}
\newtheorem{thm}{Theorem}

\theoremstyle{definition}

\theoremstyle{remark}

\begin{document}


\title{The Sum of Squares Law}

\author{Julio Kovacs\thanks{Corresponding author. Email: {\tt julio@quantumgravityresearch.org}},\hspace{.5ex} Fang Fang, Garrett Sadler, and Klee Irwin\thanks{Group leader. Email: {\tt klee@quantumgravityresearch.org}}\\[8pt]{\em Quantum Gravity Research,} Topanga, CA, U.S.}

\date{27 September 2012}

\maketitle

\pagestyle{plain} 
\setcounter{page}{1}

\begin{abstract}
We show that when projecting an edge-transitive $N$-dimensional polytope onto an $M$-dimensional subspace of $\R^N$, the sums of the squares of the original and projected edges are in the ratio $N/M$.
\end{abstract}

\section*{Statement}
Let $X\subset \R^N$ a set of points that determines an $N$-dimensional polytope. Let $E$ denote the number of its edges, and $\sigma$ the sum of the squares of the edge lengths. Let $S$ be an $M$-dimensional subspace of $\R^N$, and $\sigma'$ the sum of the squares of the lengths of the projections, onto $S$, of the edges of $X$.

Let $G$ be the group of {\em proper} symmetries of the polytope $X$ (that is, no reflections). If $G$ acts transitively on the set of edges of $X$, then:
\[
  \sigma' = \sigma\cdot\frac{M}{N}.
\]

\section*{The orthogonality relations}
The basic result used in our proof is the so-called {\em orthogonality relations} in the context of representations of groups. The form of these relations that we need is the following:

\begin{thm}
Let $\Gamma:G\ra V\times V$ be an irreducible unitary representation of a finite group $G$. Denoting by $\Gamma(R)_{nm}$ the matrix elements of the linear map $\Gamma(R)$ with respect to an orthonormal basis of $V$, we have:
\begin{equation}\label{eq:orth-rel}
  \sum_{R\in G}^{|G|} \Gamma(R)_{nm}^* \Gamma(R)_{n'm'} = \delta_{nn'}\delta_{mm'}\frac{|G|}{\dim V},
\end{equation}
\end{thm}
\noindent where the ${}^*$ denotes complex conjugation.

A proof of these relation can be found in standard books on representation theory, for instance \cite[p.\ 79]{Brocker1985} or \cite[p.\ 14]{Serre1977}. See also the Wikipedia article \url{http://en.wikipedia.org/wiki/Schur\_orthogonality\_relations}.

\section*{Proof of the sum of squares law}
The idea is apply the orthogonality relations (\ref{eq:orth-rel}) to the group $G$ of proper symmetries of the polytope $X$, considering its standard representation on the space $\R^N$ (i.e., $R\cdot x = R(x)$). This representation is clearly unitary, since the elements of the group are rotations and hence orthogonal transformations. Also, the representation is irreducible, since $G$ takes a given edge to all the other edges, which do not lie on any proper subspace due to the assumption of $X$ being an $N$-dimensional polytope.

We can assume that the edge lengths of $X$ are all equal to 1. Let $\{v_1,\dots,v_N\}$ be an orthonormal basis for $\R^N$ such that $v_1$ coincides with the direction of one of the edges $e$ of $X$. Then, for $R\in G$, let $\Gamma(R)$ be the matrix of $R$ in that basis, that is:
\[
  R(v_j) = \sum_{i=1}^N \Gamma(R)_{ij} v_i.
\]
Since this is an orthonormal basis, we have:
\[
  \Gamma(R)_{ij} = \langle R(v_j),v_i \rangle,
\]
where $\langle,\rangle$ denotes the standard inner product in $\R^N$. In particular, for $j=1$:
\begin{equation}\label{eq:sos0}
  \Gamma(R)_{i1} = \langle R(e),v_i \rangle \qquad (i=1,\dots,N).
\end{equation}
Note that this is exactly the length of the projection of each edge onto the line spanned by $v_i$. Now, from equation (\ref{eq:orth-rel}), by putting $n'=n$ and $m'=m$, we get:
\begin{equation}\label{eq:sos1}
  \sum_{R\in G} \abs{\Gamma(R)_{nm}}^2 = \frac{|G|}{N}.
\end{equation}
Using the $\Gamma$s given by the previous equation:
\begin{equation}\label{eq:sos2}
  \sum_{R\in G} \langle R(e),v_i \rangle^2 = \frac{|G|}{N}  \qquad (i=1,\dots,N).
\end{equation}

Now let $v$ be any unit vector. We'll show that the above equality holds for $v$ as it does for $v_i$. To see this, write $v$ as a linear combination of the basis vectors $v_i$: $v=sum_i a_i v_i$. Since $\norm{v}=1$, we have $\sum a_i^2=1$. Then:
\[
\begin{split}
  \sum_R \langle R(e),v \rangle^2 &= \sum_R \langle R(e),\sum_i a_i v_i \rangle^2 = \sum_R \Bigl(\sum_i a_i\langle R(e),v_i \rangle\Bigr)^2\\
  &=\sum_R\Bigl( \sum_i a_i^2 \langle R(e),v_i \rangle^2 + 2\sum_{i<j}a_i a_j \langle R(e),v_i\rangle \langle R(e),v_j\rangle\Bigr)\\
  &= \sum_i a_i^2 \sum_R \langle R(e),v_i \rangle^2 + 2 \sum_{i<j} a_i a_j \sum_R \langle R(e),v_i\rangle \langle R(e),v_j\rangle\\
  &= \frac{|G|}{N} + 2 \sum_{i<j} a_i a_j \sum_R \Gamma(R)_{i1} \Gamma(R)_{j1},
\end{split}
\]
due to eqs.~(\ref{eq:sos2}) and~(\ref{eq:sos0}). Now it turns out that the second term is 0. This is an immediate consequence of eq.~(\ref{eq:orth-rel}) with $n=i$, $n'=j$, $m=m'=1$. Therefore, the equality:
\begin{equation}\label{eq:sos3}
  \sum_R \langle R(e),v \rangle^2 = \frac{|G|}{N}
\end{equation}
holds for any unit vector $v$.

Now let $S$ be the projection subspace of dimension $M>1$, and let's denote by $P_S:\R^N\ra S$ the projection operator. Choose an orthonormal basis $\{u_1,\dots,u_M\}$ of $S$. Then:
\[
  P_S(R(e))=\sum_{i=1}^M b_i u_i,
\]
with
\[
  b_i = \langle P_S(R(e)),u_i\rangle = \langle R(e),u_i\rangle.
\]
Therefore,
\[
  \sum_{R\in G} \norm{P_S(R(e))}^2 = \sum_R \sum_i b_i^2 = \sum_R \sum_i \langle R(e),u_i\rangle^2 = \sum_{i=1}^M \Bigl( \sum_R \langle R(e),u_i\rangle^2 \Bigr) = |G|\cdot\frac{M}{N},
\]
where the last equality is because of eq.~(\ref{eq:sos3}).

To obtain the required result, we observe that $G$ can be partitioned in $E$ ``cosets'' of the same cardinality $k$, where $E$ is the number of edges of $X$. To see this, let $H=\{g\in G\ |\ g\cdot e = e\}$ be the subgroup of $G$ that leaves edge $e$ invariant. Then the coset $RH=\{g\in G\ |\ g\cdot e = R(e)\}$ is the subset of elements of $G$ that send edge $e$ to edge $R(e)$. Denote the cardinality of $H$ by $k$. Since there are $E$ edges and the action is edge-transitive, there are $E$ cosets, each of cardinality $k$. Therefore, $|G|=kE$. Denoting the edges by $e_1,\dots,e_E$, and the corresponding cosets by $C_1,\dots,C_E$ (so that $R(e)=e_l$ for $R\in C_l$), we have:
\[
\begin{split}
  \sum_{R\in G} \norm{P_S(R(e))}^2 &= \sum_{R\in \cup_{l=1}^E C_l} \norm{P_S(R(e))}^2 = \sum_{l=1}^E \sum_{R\in C_l} \norm{P_S(R(e))}^2\\
   &= \sum_{l=1}^E \sum_{R\in C_l} \norm{P_S(e_l)}^2 = \sum_{l=1}^E k\norm{P_S(e_l)}^2 = k\sum_{l=1}^E \norm{P_S(e_l)}^2.
\end{split}
\]

On the other hand, we saw that the left-hand side of this equation equals $|G|\cdot M/N$, which is $kE\cdot M/N$. Equating this to the above and canceling the factor $k$, we obtain:
\[
  \sigma' = \sum_{l=1}^E \norm{P_S(e_l)}^2 = E\cdot\frac{M}{N} = \sigma\cdot\frac{M}{N},
\]
which completes the proof.

\bibliographystyle{abbrv}
\bibliography{m-science}

\end{document}